\documentclass[11pt, reqno]{amsart}
\usepackage{amsmath, amsthm, a4, latexsym, amssymb}

\def\@rmrk#1#2{\refstepcounter
    {#1}\@ifnextchar[{\@yrmrk{#1}{#2}}{\@xrmrk{#1}{#2}}}

\makeatletter\@addtoreset{equation}{section}\makeatother

 \newfont{\bfit}{cmbxti10 scaled 2000}
 \newfont{\biggi}{cmr12 scaled 2000}

 \newcommand{\eps}{\varepsilon}

 \newcommand{\prob}{\mathbb{P}}

 \newcommand{\me}{\mathbb{E}}
 
 \renewcommand{\P}{\mathbb{P}}
 
 \newcommand{\one}{\1}

 \newcommand{\skria}{{\mathcal A}}
 \newcommand{\skrib}{{\mathcal B}}
 
 \newcommand{\skrid}{{\mathcal D}}

 \newcommand{\skrig}{{\mathcal G}}
 \newcommand{\skrih}{{\mathcal H}}

 \newcommand{\skril}{{\mathcal L}}
 \newcommand{\skrim}{{\mathcal M}}
 \newcommand{\skriw}{{\mathcal W}}

 \newcommand{\skrip}{{\mathcal P}}

 \newcommand{\skriz}{{\mathcal Z}}
 
 \newcommand{\sfrac}[2]{\mbox{$\frac{#1}{#2}$}}

\def\1{{\mathchoice {1\mskip-4mu\mathrm l}      
{1\mskip-4mu\mathrm l}
{1\mskip-4.5mu\mathrm l} {1\mskip-5mu\mathrm l}}}

\newcommand{\eq}{\begin{equation}}
\newcommand{\en}{\end{equation}}

\newenvironment{Proof}
{\vskip0.1cm\noindent{\bf Proof. }{\hspace*{0.3cm}}}%
{\nopagebreak {\hspace*{\fill}\rule{2mm}{2mm}}\\ }

{\nopagebreak {\hspace*{\fill}\rule{2mm}{2mm}}\\ }

\renewcommand{\subsection}{\secdef \subsct\sbsect}
\newcommand{\subsct}[2][default]{\refstepcounter{subsection}
\vspace{0.15cm}
{\flushleft\bf \arabic{section}.\arabic{subsection}~\bf #1  }
\nopagebreak\nopagebreak}
\newcommand{\sbsect}[1]{\vspace{0.1cm}\noindent
{\bf #1}\vspace{0.1cm}}

\newtheorem{theorem}{Theorem}[section]
\newtheorem{lemma}[theorem]{Lemma}
\newtheorem{cor}[theorem]{Corollary}

\newtheoremstyle{thm}{1.5ex}{1.5ex}{\itshape\rmfamily}{}
{\bfseries\rmfamily}{}{2ex}{}

\newtheoremstyle{rem}{1.3ex}{1.3ex}{\rmfamily}{}
{\itshape\rmfamily}{}{1.5ex}{}
\theoremstyle{rem}

\refstepcounter{subsection}

\def\thebibliography#1{\section*{References}
  \list%
  {\arabic{enumi}.}
    {\settowidth\labelwidth{[#1]}\leftmargin\labelwidth
    \advance\leftmargin\labelsep
    \parsep0pt\itemsep0pt
    \usecounter{enumi}}
    \def\newblock{\hskip .11em plus .33em minus .07em}
    \sloppy                   
    \sfcode`\.=1000\relax}

\begin{document}
\title[Large deviation  result  for the empirical locality measure of typed random geometric  graphs]
{\Large  Large deviation  result for  the empirical locality measure
of typed random geometric  graphs}

\author[Kwabena Doku-Amponsah ]{}

\maketitle
\thispagestyle{empty}
\vspace{-0.5cm}

\centerline{\sc By Kwabena Doku-Amponsah }
\renewcommand{\thefootnote}{}
\footnote{\textit{AMS Subject Classification:} 60F10, 05C80}
\footnote{\textit{Keywords: }  typed random geometric graph,
empirical locality measure, empirical degree measure,  detached
nodes.}
\renewcommand{\thefootnote}{1}

\vspace{-0.5cm}
\centerline{\textit{University of Ghana}}

\begin{quote}{\small }{\bf Abstract.}
 In  this  article for  a finite typed  random geometric graph
 we define the empirical locality
distribution, which records the number of nodes of a given type
linked to a given number of nodes of each type.  We find large
deviation principle (LDP) for the \emph{ empirical locality measure}
given the empirical pair measure and  the empirical type measure of
the typed random geometric graphs. From this LDP, we derive large
deviation principles  for the \emph{degree measure and the
proportion of detached nodes} in the classical Erd\H{o}s-R\'{e}nyi
graph defined on $[0, 1]^d.$ This graphs have been suggested by
(Canning and Penman, 2003) as a possible extension to the randomly
typed random graphs.
\end{quote}\vspace{0.3cm}

\vspace{0.3cm}
\section{Introduction}

The typed random geometric graph (TRGG)  is  obtain  when
 $n$  nodes or  nodes or points are placed uniformly at random in $[0, 1]^d,$
and in  addition each  node  is  assigned  an independently chosen
 \emph{type} or \emph{colour} or \emph{symbol} or \emph{spin} from a
finite alphabet $\skrib$.  And any two points with types $a_1,a_2$
(resp.) distance at most $r_n(a_1,a_2)$ apart are linked. The
linking radius $r_n$ plays similar role as the link probability
$p_n$ in the randomly coloured random graph models introduced by
(Penman, 1998), surveyed by (Canning and Penman, 2003) and
 studied in  (Doku-Amponsah and Moerters, 2010). The classical Erd\H{o}s-R\'{e}nyi graph on $[0, 1]^d$
 is obtain when
$n$ points are chosen at random uniformly and independently from
$[0,1]^d$ and $\lambda_n$ edges are inserted at random among the
nodes.

In  this article we  extend the LDP for the \emph{ empirical
locality measure} conditioned on the empirical pair measure  and the
empirical type distribution,  see (Doku-Amponsah et al., 2010,
Theorem~2.5), to TRGG models. From this result we derive the LDP for
the degree distribution and proportion of detached nodes the
classical Erd\H{o}s-R\'{e}nyi graph defined on $[0, 1]^d.$   See
(Doku-Amponsah, 2014)  for similar result for the classical
Erd\"{o}-Renyi graphs.

Note that the LDP for the \emph{ empirical locality measure}
conditioned on a  given empirical pair measure  and empirical type
measure  of TRGG is a  crucial  step  in  the establishment of a
full large deviation principle for the \emph{ empirical locality
measure} of TRGG. See (Doku-Amponsah, 2014[b]).

\subsection{TRGG  Model.}
The  TRGG is  a  general model of random geometric graphs  in which
the linking radius depends on the type or type of the nodes. The
main statistics for this model of random graphs are the empirical
pair distribution  and the empirical type distribution.

Given   a probability measure $\nu$ on $\skrib$ and a symmetric
function $r_n\colon\skrib\times\skrib\rightarrow (0,1]$ we  define
the {\em randomly typed random geometric graph} or  \emph{typed
random geometric graph}~$X$ having $n$ nodes as follows: Pick nodes
$X_1,...,X_n$  at  random  independently according to the uniform
distribution on $[0,\,1]^d.$ we assign to each node $X_j$ type
$Z(X_j)$ at random  and independently according to the {\em type
law} $\nu.$ Given the types, we link any two nodes
$X_i,X_j$,$(i\not=j)$ by an edge independently of everything else,
if $$\|X_i-X_j\|\le r_n\big[Z(X_i),Z(X_j)\big].$$ In this  article
we  shall  refer to $r_n(a_1,a_2),$   for $a_1,a_2\in\skrib$ as a
link radius, and always look  at
$$X=((Z(X_i),Z(X_j))\,:\,i,j=1,2,3,...,n),E)$$ under the joint measure of
graph and type. We look at  $X$ as TRGG with nodes $X_1,...,X_n$
chosen at random uniformly   and independently from the nodes space
$[0,1]^d.$ For  the  purposes of  this  study we restrict ourselves
to  the near intermediate cases .i.e. the link radius $r_n$
satisfies the condition $n r_n^d(a_1,a_2) \to \lambda(a_1,a_2)$ for
all $a_1,a_2\in \skrib$, where $\lambda\colon\skrib^2\rightarrow
[0,\infty)$ is a symmetric function, which is not identically equal
to zero.

For any  set of finite or countable elements $\skrib,$ let
$\skrip(\skrib)$ be the space of probability vectors, and
$\tilde\skrip(\skrib)$ the space of finite vectors on $\skrib$, both
equipped with the weakest topology on  $\skrip(\skrib).$  By
convention we write
$$\skriz=\{0,1,2,...\}.$$
We associate with any typed graph $X$ a probability measure, the
\emph{empirical type distribution}~$\skril_X^1\in\skrip(\skrib)$,~by
$$\skril_X^{1}(a):=\frac{1}{n}\sum_{j=1}^{n}\delta_{Z(X_j)}(a),\quad\mbox{ for $a_1\in\skrib$, }$$
and the \emph{empirical pair measure}
$\skril_X^{2}\in\tilde\skrip_*(\skrib^2),$ by
$$\skril_X^{2}(a_1,a_2):=\frac{1}{n}\sum_{(i,j)\in E}[\delta_{(X(Y_i),Z(X_j))}+
\delta_{((Z(X_j),X(Y_i))}](a_1,a_2),\quad\mbox{ for $(a,
b)\in\skrib^2$. }$$ Also we define the \emph{empirical locality
measure} $\skrim_X\in\skrip(\skrib\times\skriz)$, by
$$\skrim_X(a_1,\sigma):=\frac{1}{n}\sum_{j=1}^{n}\delta_{(Z(X_j),\skria(X_j))}(a_1,\sigma),\quad
\mbox{ for $(a_1,\sigma)\in\skrib\times\skriz$, }$$ where
$\skria(i)=(\sigma^{i}(a_2),\,a_2\in\skrib)$ and $\sigma^{i}(b)$ is
the number of nodes  of type $a_2$  linked to node $i$.

For any $n\in\skriz$ we define
$$\begin{aligned}
\skrip_n(\skrib) & := \big\{ \varpi\in \skrip(\skrib) \, : \, n\varpi(a) \in \skriz \mbox{ for all } a\in\skrib\big\},\\
\tilde \skrip_n(\skrib \times \skrib) & := \big\{ \omega\in
\tilde\skrip_*(\skrib\times\skrib) \, : \, \sfrac
n{1+\one\{a_1=a_2\}}\,\omega(a_1,a_2) \in \skriz  \mbox{ for all }
a_1,a_2\in\skrib \big\}\, ,
\end{aligned}$$

\subsection{Conditional TRGG.}\label{intro2}
 Let $\varpi(a_1)>0,$ for all $a_1\in\skrib$. We  observe that the distribution of the TRGG if the  empirical type measure $\varpi_n$ and
empirical pair distribution~$\omega_n$,
$$\prob_{(\varpi_n,\omega_n)}:=\prob\{ \,\cdot\,  \,|\,\psi(\skrim_X)=(\varpi_n,\omega_n)\},$$
may be obtained as follows:
\begin{itemize}
\item Pick nodes  $X_1,...,X_n$  at  random
independently according to the uniform distribution on $[0,\,1]^2.$
\item Give types to the nodes by picking without replacement from the
constellation  of $n$~types, containing each type $a_1\in\skrib$
precisely $n\varpi_n(a)$ times;
\item For each  pair $\{a_1,a_2\}$ of types make precisely $n(a_1,a_2)$ links by
picking  without replacement from the pool of potential edges
linking nodes of type $a_1$ and $a_2$,
where $$m_n(a_1,a_2):=\left\{ \begin{array}{ll} n\, \omega_n(a_1,a_2) & \mbox{if }a_1\not=a_2\\
\frac n2\, \omega_n(a_1,a_2) & \mbox{if }a_1=a_2\,
.\end{array}\right.$$
\end{itemize}

In  the remainder of the paper we state and prove  our  LDP results.
 In Section~\ref{mainresults}  we   state our LDPs,
 Theorem~\ref{maing3},  Theorem~\ref{maing4}, and Corollary~\ref{maing5}.
  The proof of Theorem~\ref{maing3}, carried out in
Section~\ref{proofmaing3}, uses a combinatorial arguments based on
random allocation of typed balls into typed bins  and
(Doku-Amponsah,Lemma~5, 2014). The article ends with the proofs of
Theorem~\ref{maing4}~and~\ref{maing5}  in
Subsection~\ref{proofcorollaries}.

\section{Statement of the results}\label{mainresults}
For any $\ell\in\skrip(\skrib\times\skriz^{\skrib})$we  denote by
$\ell_1$ the $\skrib-$ marginal of $\ell$ and for  every $(a_2,
a_1)\in\skrib\times\skrib,$ let $\ell_2$  be the distribution   of
the couple $(a_1,\sigma(a_2))$ under  the  measure $\ell.$ We define
the finite  measure, $\skrih\in\tilde\skrip(\skrib\times\skrib)$ by
$$\skrih_2(\ell)(a_2, a_1):=
\sum_{\sigma(a_2)\in\skriz}\ell_2(a_1,\sigma(a_2))\sigma(a_2),
\quad\mbox{ for $a_1,a_2\in\skrib$}$$ and  write
$\skrih_1(\ell)=\ell_1.$ We define the function  $\skrih \colon
\skrip((\skrib\times\skriz^{\skrib}) \to \skrip(\skrib) \times
\tilde\skrip(\skrib \times \skrib)$ by
$\skrih(\ell)=(\skrih_1(\ell),\skrih_2(\ell))$ and observe that
$\skrih(\skrim_X)=(\skril_X^1, \skril_X^2).$ Note  that, in
 the  weak  topology $\skrih_1$ is not discontinuous function  but $\skrih_2$
 is. To  be specific, in  the expression $\displaystyle
\sum_{\sigma(a_2)\in\skriz}\ell_2(a_1,\sigma(a_2))\sigma(a_2)$ the
function $\sigma(a_2)$ may not be bounded and hence in  the  weak
topology the functional $\displaystyle \ell\to\skrih_2(\ell)$ would
be discontinuous. We say two of measures
$(\omega,\ell)\in\tilde{\skrip}(\skrib\times\skrib)\times\skrip(\skrib\times\skriz^{\skrib})$
\emph{consistent} if
\begin{equation}\label{consistent}
\skrih_2(\ell)(a_2, a_1) = \omega(a_2, a_1), \quad\mbox{ for all
$a_1,a_2\in\skrib.$}
\end{equation}

 The next theorem gives LDP for the empirical locality measure of a sequence of graphs with
given empirical type  distribution and empirical pair distribution.

\begin{theorem}\label{randomg.LDprobg}\label{maing3} Let the sequence
$(\varpi_n,\omega_n)\in\skrip_n(\skrib)\times\tilde{\skrip}_{n}(\skrib\times\skrib)$
converges to a limit
$((\varpi,\omega)\in\skrip(\skrib)\times\tilde{\skrip}_*(\skrib\times\skrib).$
Suppose that $X$ is a TRGG  graph  conditioned on the  set $\{
\skrih(\skrim_X)=(\varpi_n,\omega_n) \}$. Then, the empirical
locality measure  $\skrim_X$, as $n\rightarrow\infty,$ satisfies  an
LDP  in the space $\skrip(\skrib\times\skriz)$ with good rate
function
\begin{align}\label{randomg.rateLDprob}
\tilde{J}_{((\varpi,\omega)}(\ell)=\left\{
\begin{array}{ll}H(\ell\,\|\,Q_{poi}) & \mbox
  {if  $(\omega,\ell)$  is  consistent and $\ell_1=\varpi_2$ }\\
\infty & \mbox{otherwise.}
\end{array}\right.
\end{align}

where  $$Q_{poi}(a_1\,,\,\sigma)=\ell_{1}(a_1)\prod_{a_2\in\skrib}
e^{-\frac{\omega(a_1,a_2)}{\ell_1(a_1)}} \, \frac{1}{\sigma(a_2)!}\,
\Big(\frac{\omega(a_1,a_2)}{\ell_1(a)}\Big)^{\sigma(a_2)},
\quad\mbox{for $a_1\in\skrib$, $\sigma\in\skriz$} .$$

\end{theorem}

Note that \emph{degree distribution} $D_X \in \skrip(\skriz)$ of a
graph with empirical locality distribution $\skrim_X$ is given   by
$$D_X(r)= \sum_{a_1\in\skrib}\sum_ {\sigma\in\skriz} \delta_r\big({\textstyle \sum_{a_2\in\skrib}}
 \sigma(a_2)\big)\,
\skrim_X(a_1,\sigma), \qquad \mbox{ for $r\in\skriz$,}$$ i.e.
$D_X(r)$ is the proportion of nodes in the graph with degree~$r$.
Theorem~\ref{maing4}  below  is  a  spacial  case  of
Theorem~\ref{maing3}  above  where  $\skrim_X=D_X$ ,  the degree
distribution  and $\langle \skrih(\skrim_X)\rangle =2|E|/n.$

We  write
$$\rho(d)=\sfrac{\pi^{d/2}}{\Gamma\big(\sfrac{(d+2)}{2}\big)},$$
  where  $\Gamma$  is  the  gamma  function.

\begin{theorem}\label{randomg.LDprobg}\label{maing4} Suppose the sequence
$\lambda_n/n$ converges to a limit $\rho(d)t/2.$ Let
$\skrig(n,\lambda_n)$ be a Random  geometric graph, where nodes
$X_1,...,X_n$ are chosen at random uniformly and independently from
$[0,1]^d,$ and $\lambda_n$ edges  are  inserted  at random among the
nodes. Then, as $n\rightarrow\infty,$ the degree distribution $D_X$
of $\skrig(n,\lambda_n)$ satisfies  large  deviation principle on
the space $\skrip(\skriz)$ with good rate function
\begin{align}\label{randomg.rateLDprob}
\eta(\delta)=\left\{
\begin{array}{ll}H(\delta\,\|\,q_{\langle \delta\rangle }) ,\,& \mbox{$\langle \delta\rangle =\rho(d)t$ }\\
\infty & ,\mbox{otherwise.}
\end{array}\right.
\end{align}
\end{theorem}

From  Theorem~\ref{maing4} above  we  obtain the  following
Corollary~\ref{main6}  using  the  contraction principle.See,
\cite{DZ98}.

\begin{cor}\label{main6}\label{maing5}
Suppose the sequence $\lambda_n/n$ converges to a limit
$\rho(d)t/2.$ Let $\skrig(n,\lambda_n)$ be a random geometric graph,
where nodes $X_1,...,X_n$ are chosen at random uniformly and
independently from $[0,1]^d,$ and $\lambda_n$ edges  are inserted at
random among the nodes. Then, the proportion of detached nodes
$D_X(0)$ of $\skrig(n,\lambda_n),$  as $n\to\infty,$ obeys an LDP on
the  space $[0,1]$ with good, convex rate function
\begin{align}\label{randomg.detached1}
\xi(y)= \left\{ \begin{array}{ll}
y\log\sfrac{y}{1-e^{-\rho(d)t}}+(1-y)\log\sfrac{(1-y)}{(1-e^{-\rho(d)t})}+\rho(d)t\log\lambda-\rho(d)t\log
c\rho(d), & \mbox{if\,
$y\ge 1-t\rho(d) $,}\\
\infty & \mbox{ if  $ y< 1-t\rho(d) $,}
\end{array}\right.
\end{align}
where  $\alpha=\alpha(y,t)$ uniquely  solve the  equation
$\sfrac{1-e^{-\alpha}}{\alpha}=\sfrac{1-y}{\rho(d)t}.$
\end{cor}

(Doku-Amponsah,2014) obtained the same  result for the proportion of
detached nodes in $\skrig(n,nt/2).$

\section{Proof of the Results}\label{proofmaing3}
\subsection{Proof of  Theorem~\ref{maing3} by  Random Allocation.}
We  recall  the  conditional TRGG from  Subsection~\ref{intro2} and
 denote by $\skriw(a_1)$ the collection of nodes  which have type
$a\in\skrib$. Note that $$\sharp \skriw(a_1)=n\varpi_n(a).$$

The random allocation model is appropriately  obtained  when  typed
balls are dropped at random in typed bins.
 In the next Lemma we prove the exponential equivalence for the
 distribution
of ~$\tilde{\skrim}_X$, see (Doku-Amponsah, 2014), with respect to
$\tilde{\prob}_{(\varpi_n,\omega_n)}$ the law of the random
allocation model and $\skrim$ with respect to
$\prob_{(\varpi_n,\omega_n)}$  where
$$\prob_{(\varpi_n,\omega_n)}=\prob\{ \,\cdot\,  \,|\,\skrih(\skrim_X)=(\varpi_n,\omega_n)\}.$$
Thus the distribution of the TRGG conditioned to have type law
$\varpi_n$ and edge distribution $\varpi_n$. Recall the definition
of exponential equivalence, see (Dembo et al.,1998,
Definition~4.2.10).We define the metric $\skrid$ of total variation
by
$$\skrid(\ell,\tilde{\ell})=\sfrac{1}{2}\sum_{(a_1,\sigma)\in\skrib\times\skriz^{\skrib}}
|\ell(a_1,\sigma)-\tilde{\ell}(a_1,\sigma)|, \quad \mbox{ for
}\ell,\tilde{\ell}\in\skrip(\skrib\times\skriz^{\skrib})$$  and
observe that this metric generates the weak topology.

\begin{lemma}\label{randomg.expequivalnce}
For every $\eps>0,$
\begin{equation}\label{randomg.totalv}
\lim_{n\rightarrow\infty}\sfrac{1}{n}\log\prob\big\{
\skrid(\skrim_X,\tilde{\skrim}_X)\ge\eps\big\}=-\infty,
\end{equation}
where $\prob$ is  a suitable coupling  between  the  laws
$\tilde{\prob}_{(\varpi_n,\omega_n)}$ and
$\prob_{(\varpi_n,\omega_n)}.$
\end{lemma}

\begin{Proof} Proof  of  this  Lemma  given  below is  also uses the same
coupling argument  of (Boucheron et al., 2003) presented  in
(Doku-Amponsah, 2014).

For each $a_1,a_2\in\skrib$ we begin as  follows: At every step
$k=1,\ldots, m_n(a_1,a_2),$ we pick  at  random two nodes
$\skriw^k_i\in \skriw(a_1)$ and $\skriw^k_j\in \skriw(a_2)$.  we
place one ball of type $a_2$ into bin $i$ and one ball of type $a_1$
in $j,$ and join $\skriw^k_i$ with $\skriw^k_j$ by an edge except
 when $\skriw^k_i=\skriw^k_j$ or the pair of nodes already formed an  edge. If one
of these two scenario occur, then we only pick a link at  random
from the set of all likely edges joining types $a_1$ and $a_2$,
which are not already an  edge in our graph model. This ends the
formation of the graph having  $\skril_X^1=\varpi_n$,
$\skril_X^2=\omega_n$. For this collection  denote, for each bin
$i\in\{1,\ldots, n\}$, by $Z(X_i)$ its type, and by $\sigma^i(a_2)$
the number of balls of type $a_2\in\skrib$ it contains, and  define
the \emph{empirical occupancy measure} of the collection by
$$\tilde{\skrim}_X(a_1,\sigma)
= \frac 1n \sum_{v\in V} \delta_{(\tilde{X}(v),
\tilde{\skria}(v))}(a_1, \sigma), \qquad \mbox{ for }
(a_1,\sigma)\in\skrib\times\skriz^{\skrib}.$$

Observe that
\begin{equation}\label{randomg.XTX}
\skrid(\skrim_X\,,\,\tilde{\skrim}_X)\le \sfrac{2}{n}
\sum_{a_1,a_2\in\skrib}B^{n}(a_1,a_2)\, ,
\end{equation}
where $B^n(a_1,a_2)$ is the total number of steps $u\in\{1,\ldots,
m_n(a_1,a_2)\}$ at which there is inconsistency between the nodes
$\skriw^u_i$, $\skriw^u_j$ sampled and the nodes that received the
$u^{\rm th}$ edge linking $a_1$ and $a_2$ in the graph formation.

If $a_1,a_2\in\skrib$,the frequency of $\skriw^k_i=\skriw^k_j$ or
the two nodes are already linked is equal to
$$p_{[k]}(a_1,a_2):=\sfrac{1}{m_n(a_1,a_2)}\1_{\{a_1=a_2\}}+\big(1-\sfrac{1}{m_n(a_1,a_2)}
\1_{\{a_1=a_2\}}\big)\sfrac{(k-1)}{(m_n(a_1,a_2))^2}.$$
$B^{n}(a_1,a_2)$ is a sum of independent $0~or~1$ random variables
$X_1,\,...,\,X_{m_n(a_1,a_2)}$ with `success' frequencies  equal to
$p_{[1]}(a_1,a_2), \ldots, p_{[m_n(a_1,a_2)]}(a_1,a_2)$. Note that
$\me[X_k]= p_{[k]}(a_1,a_2)$  and
$$Var[X_k]=p_{[k]}(a_1,a_2)(1-p_{[k]}(a_1,a_2)).$$  Now,  we   have   $$\me B^{n}(a_1,a_2)= \sum_{k=1}^{m_n(a_1,a_2)} p_{[k]}(a_1,a_2)=\1_{\{a_1=a_2\}}
+\big(1-\1_{\{a_1=a_2\}}\sfrac{1}{m_n(a_1,a_2)}\big)\big(1-\sfrac{1}{m_n(a_1,a_2)}\big)\le
1+\1_{\{a_1=a_2\}}.$$

We  write
$$\sigma_n^2(a_1,a_2):=\sfrac{1}{m_n(a_1,a_2)}\sum_{k=1}^{m_n(a_1,a_2)}Var[X_k]$$
and  observe  that   $$\lim_{n\to
\infty}\me(B^n(a_1,a_2))=\lim_{n\to
\infty}Var(B^n(a_1,a_2))=\lim_{n\to
\infty}m_n(a_1,a_2)\sigma_n^2(a_1,a_2)=\1_{\{a_1=a_2\}}+1.$$

Define  $h(t)=(1+t)\log(1+t)-t,$ for  $t\ge 0$   and apply Bennett's
inequality, see ( Bennett,  2004)  to  obtain, for very large $n$
$$\P\big\{ \sfrac{ B^{n}(a_1,a_2)}{n}\ge\sfrac{ \1_{\{a_1=a_2\}}+1}{n}+\delta_{1}\big\}
\le
exp\Big[-m_n(a_1,a_2)\sigma_n^2(a_1,a_2)h(\sfrac{n\delta_{1}}{n(a_1,a_2)\sigma_n^2(a_1,a_2)})\Big],$$
for any $\delta_1>0.$ Let $\eps\ge 0 $ and  choose
$\delta_1=\sfrac{\eps}{2m^2}.$ Suppose that we have
$B^n(a_1,a_2)\le\delta$. Then, by~\eqref{randomg.XTX},
$$d(\skrim_X,{\ell}_n)\le 2\delta_1 m^2=\eps.$$ Hence,
$$\begin{aligned}
\prob\big\{ \skrid(\skrim_X,\tilde{\skrim}_X) > \eps \big\}  \le
\sum_{a_1,a_2\in\skrib} \prob\big\{ B^n(a_1,a_2)\ge n\delta_1
\big\}&\le m^2\sup_{a_1,a_2\in\skrib}\prob\big\{ B^n(a_1,a_2)\ge
\1_{\{a_1=a_2\}}+1+ (n\delta_1)/2
\big\}\\
 & \le m^2\sup_{a_1,a_2\in\skrib} exp\Big[-m_n(a_1,a_2)\sigma_n^2(a_1,a_2)h(\sfrac{n\delta_{1}}{m_n(a_1,a_2)\sigma_n^2(a_1,a_2)})\Big] .
\end{aligned}$$

Let   $0\le \delta_2\le 1$. The,  for  very  large $n$ we that
 have
\begin{equation}\begin{aligned}\label{Equ.coupling}
\frac{1}{n}&
\log\P\Big\{\skrid(\skrim_X,\tilde{\skrim}_X) > \eps \Big\}\le-(1-\delta_{2})h(\sfrac{n\delta_1}{2(1+\delta_{2})})\\
&=-(\1_{\{a_1=a_2\}}+1-\delta_{2})\Big[(\sfrac{1}{n}+\sfrac{\delta_1}{2(\1_{\{a_1=a_2\}}+1+\delta_{2})})\log(1+\sfrac{n\delta_1}{2(\1_{\{a_1=a_2\}}+1+\delta_{2})})-\sfrac{\delta_1}{2(\1_{\{a_1=a_2\}}+1+\delta_{2})}\Big].
\end{aligned}\end{equation}
This ends the proof of the Lemma.
\end{Proof}

To  conclude  the  proof  of Theorem~\ref{maing3},  we  note  that
empirical occupancy  measure  $\tilde{\skrim}_X$  is  exponential
equivalent  to  $\skrim_X,$ and further    $\tilde{\skrim}_X$ under
the law $\tilde{\prob}_{(\varpi_n,\omega_n)}$  obeys a  large
deviation principle with  rate function
$\tilde{J}_{((\varpi,\omega)}$ by (Doku-Amponsah, Lemma~5, 2014).
Therefore,  by  the exponential equivalent theorem, see (Dembo et
al.,1998, Theorem~4.2.13), $M$ obeys an LDP with the rate function
$\tilde{J}_{((\varpi,\omega)}.$


 \subsection{Proof of Corollary~\ref{main6} by the  Contraction  Principle.}\label{proofcorollaries}We  prove
 Corollary~\ref{main6}  from   Theorem~\ref{maing4}  by using  the contraction
 principle, \cite[Theorem~4.2.1]{DZ98} on the
linear mapping $G:\skrip(\skriz)\to [0,1]$ given  by
$G(\delta)=\delta(0).$ To  be  specific,  Theorem~\ref{maing4}
implies an LDP for random variable $G(D_X)=D_X(0)$ with  good  rate
function
$$\xi(y)=\inf\big\{H(\delta\,\|\,q_{\langle \delta\rangle} ):\delta\in\skrip(\skriz),
\delta(0)=y, \,\sum_{r=0}^{\infty}r\delta(r)= \rho(d)t\,\big\}.$$

Note that,  we have
$$\rho(d)t=\sum_{k=1}^{\infty}k\delta(k)\ge
\sum_{k=1}^{\infty}\delta(k)=1-y,$$  and  that the  class of
measures  satisfying the  two constraints is necessarily empty if
$t\rho(d)<1-y$. If $\rho(d)t\ge 1-y$,using the Lagrangian method we
can  calculate  the minimizer  $p,$ defined by $ p(0)=y$,
$p(k):=U(y,t)^{-1}\sfrac{(\alpha(y,t))^k}{k!}$ where $\alpha(y,c)$
uniquely  solve
$$\sfrac{e^{\lambda}-1}{\lambda}=\sfrac{1-y}{\rho(d)t}$$
and  $U(y,t):=\sfrac{e^{\lambda-1}}{1-y}.$  Therefore, we  have

\begin{equation}\begin{aligned}\label{equ.last}
\eta(y)&=x\log\frac{y}{q_{c\rho(d)}(0)}
+(1-y)\log\frac{(1-y)}{1-q_{c\rho(d)}(0)}+(1-y)\sum_{k=1}^{\infty}p(k)\log\sfrac{p(k)}{\hat{q}_{\rho(d)t}(k)}.\\
&=y\log\frac{y}{q_{\rho(d)t}(0)}
+(1-y)\log\frac{(1-y)}{1-q_{t\rho(d)}(0)}+t\rho(d)\log
\sfrac{\lambda}{t\rho(d)}\\
&=y\log\frac{y}{e^{-\rho(d)t}}+(1-y)\log\frac{(1-y)}{1-e^{-\rho(d)t}}+t\rho(d)\log
\sfrac{\lambda}{t\rho(d)}
\end{aligned}
\end{equation}
if  $t\rho(d)\ge 1-y$  and  $\infty$  otherwise. Particularly, if
$y=e^{-t\rho(d)}$ then we  have $\alpha(y,t)=t\rho(d),$  which gives
$\eta(e^{-t\rho(d)})=0.$ This  end  the  proof  of
Theorem~\ref{main6}.

\bigskip


\end{document}